\documentclass[12pt]{amsart}
\usepackage{amssymb}
\usepackage{amsmath}
\numberwithin{equation}{subsection}
\def\R{\mathbb R}
\def\C{\mathbb C}

\def\D{\mathbb D}
\def\CC{\widehat{\mathbb C}}
\def\N{\mathbb N}
\def\FF{{\mathcal F}}

\def\OO{{\mathcal O}}

\def\re{\operatorname{Re}}

\newtheorem{thm}[equation]{}

\theoremstyle{remark}
\newtheorem{exam}[equation]{}

\newtheorem*{rem}{Remark}
\newtheorem*{acknowledgement}{Acknowledgment}
\begin{document}
\title{Bloch's principle}
\author{Walter Bergweiler}
\email{bergweiler@math.uni-kiel.de}
\thanks{Supported by the G.I.F.,
the German--Israeli Foundation for Scientific Research and
Development, Grant G -809-234.6/2003, 
and by the Alexander von Humboldt Foundation.}
\address{Mathematisches Seminar,
Christian--Albrechts--Universit\"at zu Kiel,
Ludewig--Meyn--Str.\ 4,
D--24098 Kiel,
Germany}
\begin{abstract}
A heuristic principle attributed to Andr\'e Bloch says that 
a family of holomorphic functions is likely to be normal if 
there is no nonconstant entire functions with this property.
We discuss this principle and survey recent results that have been
obtained in connection with it.
We pay special attention to properties related
to exceptional values of derivatives and existence of 
fixed points and periodic points,  but we also discuss some
other instances of the principle.
\end{abstract}
\keywords{Normal family, quasinormal, Zalcman lemma, exceptional value,
fixed point, periodic point}
\subjclass{Primary 30D45; Secondary 30D20, 30D35}
\maketitle
\tableofcontents
\section{Bloch's  heuristic principle}\label{intro}
\subsection{Introduction}\label{intro1}
A family of meromorphic functions is called normal if every sequence in the family
has a subsequence which converges 
(locally uniformly with respect to the spherical metric).
The concept of a normal family was introduced
already in 1907 by
P.~Montel~\cite{Mon07}, but there has been a lot of interest in normal families again in recent years,  
an important factor being their central role in complex dynamics.

One guiding principle in their study 
has been the heuristic principle which 
says that a family of functions meromorphic (or holomorphic) in a domain and possessing a certain 
property is likely to be normal if there is no nonconstant function meromorphic (holomorphic)
in the plane which has this property. 
This heuristic principle is usually attributed to A.~Bloch,
but it does not seem to have been stated explicitly by him
-- although his statement
``Nihil est in infinito quod non prius fuerit in finito'' made 
in his 1926 papers~\cite[p.~84]{Blo26a} and~\cite[p.~311]{Blo26b} may be 
interpreted this way.
The first explicit formulation of the heuristic principle
seems to be due to G.~Valiron~\cite[p.~2]{Val29} in 1929.
However, in~\cite[p.~4]{Val37} Valiron mentions Bloch in this context.

Here we survey some of the results that have been 
obtained in connection with Bloch's principle.
We concentrate on properties related to exceptional values of derivatives
in~\S 2 and on properties related to fixed points of iterates in~\S 3.
Excellent references for Bloch's principle are~\cite[Chapter~4]{Schi} and~\cite{Zal98}.

We mention that Andr\'e Bloch's life was tragic. He murdered his brother, uncle and aunt in 1917.
Judged not responsible for his actions, he was confined to a psychiatric hospital where
he remained until his death in 1948. It was there that he did all
his mathematical work,
including that cited above. For a detailed account of Bloch's life and work we refer
to~\cite{CarFer,Val49}.

\subsection{The basic examples}\label{basic}
A simple example for Bloch's principle is given by  the property of being bounded.
Liouville's Theorem says that a bounded entire function is constant;
that is, if $f:\C\to\C$ is holomorphic and if there exists a constant
$K$ such that $|f(z)|\leq K$ for all $z\in\C$, then $f$ is constant.
And in fact a family $\FF$ of  functions holomorphic in a domain $D$ is normal 
if there exists a constant $K$ such that $|f(z)|\leq K$ for all $z\in D$ and all $f\in\FF$.
Note, however, that the family $\FF$ of  all functions $f$ holomorphic in a domain $D$ 
for which there exists a constant $K=K(f)$ such that $|f(z)|\leq K(f)$ for all $z\in D$ need not 
be normal.
In fact, let $D$ be the unit disk $\D$, that is,
$D=\D:=\{z\in\C:|z|< 1\}$, and let $\FF:=\{f_n:n\in\N\}$,
where $f_n(z):=nz$. Each $f_n$ is bounded in $D$ (by the constant~$n$), but it is easily
seen that $\FF$ is not normal.
So in order that the property ``$f$ is bounded'' leads to a normal family it is essential
that the bound not depend on the function $f$. This example already
shows that some care has to
be taken when formulating Bloch's principle.

A key 
example deals with  the property 
of omitting three points in the Riemann sphere
$\CC:=\C\cup\{\infty\}$. 
Perhaps the most important theorem in the whole subject of normal families
is Montel's result that this condition implies normality.
\begin{thm}{\bf Montel's Theorem.}\label{thmon}
Let $a_1,a_2,a_3\in \CC$ be distinct, let $D\subset \C$ be a domain
and let $\FF$ be a family of functions meromorphic
in $D$ such that $f(z)\neq a_j$ for all $j\in\{1,2,3\}$, all $f\in\FF$, and 
all $z\in D$. Then $\FF$ is normal.
\end{thm}
The analogous statement about functions in the plane is Picard's theorem.
\begin{thm}
{\bf Picard's Theorem.}\label{thpic}
Let $a_1,a_2,a_3\in \CC$ be distinct and let $f:\C\to\CC$ be meromorphic.
If $f(z)\neq a_j$ for all $j\in\{1,2,3\}$ and 
all $z\in \C$, then $f$ is constant.
\end{thm}

In Montel's Theorem it is again essential that the points $a_1,a_2,a_3$ not depend 
on~$f$. We mention, however, that the condition
$f(z)\neq a_1(f),a_2(f),a_3(f)$ still implies normality if there exists $\varepsilon >0$
such that $\chi(a_j(f),a_k(f))\geq \varepsilon $ for $j\neq k$ and all $f$, where 
$\chi(\cdot,\cdot)$ denotes the spherical distance; see, e.~g.,~\cite[p.~104]{Schi}
or~\cite[p.~224]{Zal98}.
\subsection{The theorems of Arzel\`a-Ascoli and Marty}\label{arzela}
One basic result in the theory of normal families is a classical theorem due to
 Arzel\`a and Ascoli~\cite[p.~35]{Schi}, which we record in the following form.
\begin{thm}
{\bf Arzel\`a-Ascoli Theorem.}
A family of meromorphic functions is normal if and only if it 
is locally equicontinuous
(with respect to the spherical metric).
\end{thm}
A consequence of this is the following theorem due to Marty~\cite[p.~75]{Schi}.
Denote by 
$$f^\#(z):=\frac{|f'(z)|}{1+|f(z)|^2}=\frac12 \lim_{\zeta\to z}\frac{\chi(f(\zeta),f(z))}{|\zeta-z|}$$
the spherical derivative of a meromorphic function~$f$.
\begin{thm}
{\bf Marty's Criterion.}
A family $\FF$ of functions meromorphic in a domain $D$ is normal if and only if 
the family $\{f^\#:f\in\FF\}$ is locally bounded; that is, if for every $z\in D$ 
there exists a neighborhood $U$ of $z$ and a constant $M$ such that
$f^\#(z)\leq M$ for all $z\in U$ and for all $f\in\FF$.
\end{thm}

\subsection{Zalcman's Lemma}\label{seczl}

In order to turn the heuristic principle for certain properties into a rigorous theorem,
L.~Zalcman~\cite{Zal75} proved the following result.
\begin{thm} \label{lemzalc}
{\bf Zalcman's Lemma.}
Let $D\subset\C$ be a domain and let $\FF$ be a family of functions
meromorphic in $D$. If $\FF$ is not normal, then there exist a 
sequence $(z_k)$ in $D$, a 
sequence $(\rho_k)$ of positive real numbers, a 
sequence $(f_k)$ in $\FF$,
a point $z_0\in D$ and a nonconstant meromorphic function
$f:\C\to\CC$ 
such that
$z_k\to z_0$,
$\rho_k\to 0$ and 
$f_k(z_k+\rho_k z)\to f(z)$ locally uniformly in~$\C$.
Moreover,  
$f^\#(z)\leq f^\#(0)=1$ for all $z\in \C$.
\end{thm}
The corresponding result for normal {\em functions} rather than normal {\em families}
had been proved earlier by A.~J.~Lohwater and
Ch.~Pommerenke~\cite{LP}.
We remark that the statement 
about the spherical derivative
does not appear in~\cite{Zal75}, but it follows immediately from the proof;
see also~\cite[p.~216f]{Zal98}. This observation
plays an important role in some of the 
more recent applications; see \S\ref{seczhp} below. 

{\em To prove Zalcman's Lemma},
suppose that $\FF$ is not normal.
By Marty's Criterion, there exists a sequence $(\zeta_k)$
in $D$  tending to a point $\zeta_0\in D$ 
and a sequence $(f_k)$ in $\FF$ such that
$f_k^\#(\zeta_k)\to\infty$.
Without loss of generality, we may assume that $\zeta_0=0$ and that 
$\{z:|z|\leq 1\}\subset D$. 
Choose $z_k$ satisfying $|z_k|\leq 1$ such that 
$M_k:=f_k^\#(z_k)(1-|z_k|)=\max_{|z|\leq 1}f_k^\#(z)(1-|z|)$.
Then $M_k\geq f_k^\#(\zeta_k)(1-|\zeta_k|)$ and hence
$M_k\to\infty$.
Define $\rho_k:=1/f_k^\#(z_k)$.
Then $\rho_k\leq 1/M_k$ so that $\rho_k\to 0$.
Since $|z_k+\rho_k z|<1$ for $|z|<(1-|z_k|)/\rho_k=M_k$ 
the  function $g_k(z):=f_k(z_k+\rho_k z)$ is defined for
$|z|<M_k$ and satisfies 
$$g_k^\#(z)=\rho_k f_k^\#(z_k+\rho_k z)
\leq
\frac{1-|z_k|}{1-|z_k+\rho_k z|} 
\leq 
\frac{1-|z_k|}{1-|z_k|-\rho_k |z|} 
=\frac{1}{\displaystyle 1-\frac{|z|}{M_k}}$$
there.
By Marty's Criterion, the sequence $(g_k)$  is normal in $\C$
and thus has a subsequence
which converges locally uniformly in $\C$. Without loss of 
generality, we may assume that $g_k\to f$ for some $f:\C\to\CC$
and $z_k\to z_0$ for some $z_0\in D$.
Since $g_k^\#(0)=1$ for all $k$, we have $f^\#(0)=1$, so that
$f$ is nonconstant.
Clearly, we also have $f^\#(z)\leq 1$ for all $z\in\C$.\qed

\begin{rem}
It is not difficult to see that if $z_k,\rho_k,f_k$ and $f$ are as in 
Zalcman's Lemma, then $\FF$ is not normal at $z_0=\lim_{k\to\infty}z_k$.
In turn, we may achieve $z_k\to z_0$ for any point $z_0\in D$ at which 
$\FF$ is not normal.

The important point in 
Zalcman's Lemma is that the limit function $f$ is nonconstant.
Regardless of whether $\FF$ is normal or not 
it is always possible to choose $z_k,\rho_k$ and $f_k$ such that 
$f_k(z_k+\rho_k z)$ tends to a constant.
\end{rem}
\subsection{Zalcman's formalization of the heuristic principle}\label{seczhp}
We first note how Zalcman's lemma can be used to deduce Montel's Theorem
from Picard's Theorem. In fact, suppose that $\FF$ is a non-normal family of 
meromorphic functions which omit three values $a_1,a_2,a_3$. Let
$f_k,z_k,\rho_k$ and $f$ be as in Zalcman's Lemma. Then
obviously $f_k(z_k+\rho_k z)\neq a_j$ for all $z,k$ and $j$. Thus 
$f(z)\neq a_j$ for all $z\in\C$ by Hurwitz's Theorem, contradicting 
Picard's Theorem.

To see for what kind of properties this kind of argument may be used
we first need
to specify what a property is. 
To this end we simply collect all functions enjoying a certain property in a 
set.
It turns out to be useful to display a function $f:D\to\CC$ together with its domain
$D$ of definition. Following
A.~Robinson~\cite[\S 8]{Ro} and L.~Zalcman~\cite{Zal75}
we thus write $\langle f,D\rangle \in P$ if $f$ has the property $P$ 
on a domain 
$D$. Bloch's principle then asserts that the following two statements
should be equivalent:
\begin{itemize}
\item[$(a)$] {\em If $\langle f,\C\rangle \in P$, then $f$ is constant.}
\item[$(b)$] {\em The family $\{f: \langle f,D\rangle \in P\}$ is normal on $D$ for 
each domain $D\subset \C$.}
\end{itemize}
We say that $P$ is a {\em Bloch property} if $(a)$ and $(b)$ are equivalent.
Of course, that two statements are equivalent simply means that either
both are true or both are false. But we will later meet properties $P$ where
we can prove that $(a)$ and $(b)$ are equivalent, but where we do not know whether
either is true or false!
If statements $(a)$ and $(b)$ are actually true, then we say that $P$ is a
{\em Picard-Montel} property.

Zalcman's Lemma now leads to the following result.
\begin{thm} \label{zalcprinc}
{\bf Zalcman's Principle.}
Suppose that a property $P$ of meromorphic functions satisfies 
the following three conditions:
\begin{itemize}
\item[$(i)$] If $\langle f,D\rangle \in P$, then $\langle f|_{D'},D'\rangle \in P$ 
for every domain $D'\subset D$. 
\item[$(ii)$] If $\langle f,D\rangle \in P$ and $\varphi(z)=\rho z +c$ 
with $\rho,c\in\C$, $\rho\neq 0$,
then $\langle f\circ \varphi,\varphi^{-1}(D)\rangle \in P$.
\item[$(iii)$] Suppose that $\langle f_n,D_n\rangle \in P$ for $n\in\N$,
where $D_1\subset D_2\subset D_3\subset\dots $
and $\bigcup_{n=1}^\infty D_n=\C$.
If $f_n\to f:\C\to\CC$ locally uniformly in $\C$, then $\langle f,\C\rangle \in P$.
\end{itemize}
Then $P$ is a Bloch property.
\end{thm}
The proof that $(a)$ implies $(b)$ is a straightforward
application of the Zalcman
Lemma. In order to see that $(b)$ implies $(a)$ it suffices to consider the family 
$\{f(nz): n\in\N\}$ and note that his family is not normal if $f$ is 
meromorphic in the plane and not constant. In fact we see that conditions $(i)$ and $(ii)$
suffice in order to deduce $(a)$ from $(b)$.

The argument shows that $(b)$ follows not only from condition $(a)$, but in fact 
from the following weaker condition
\begin{itemize}
\item[$(c)$] {\em If $\langle f,\C\rangle \in P$ and if $f$ has bounded spherical derivative, 
then $f$ is constant.}
\end{itemize}
In particular, for properties $P$ satisfying the hypothesis of Zalcman's Principle
we see that $(a)$ and $(c)$ are equivalent. So in order to prove a result for functions 
meromorphic in the plane, it sometimes suffices to prove it for functions with bounded
spherical derivative. This kind of argument is due to X.~Pang,
and appears in writing first in~\cite{BE,CF,Zal94}, 
see Theorem~\ref{thm6} in \S\ref{diffpol} below. This argument 
will also appear in~\S\ref{other}, \S \ref{genzal} and~\S \ref{multiple} below.

Already here we note that the Ahlfors-Shimizu form of the Nevanlinna characteristic 
$T(r,f)$ shows
that if $f$ is meromorphic in the plane and has bounded spherical derivative, then
$T(r,f)=\OO(r^2)$ as $r\to\infty$. In particular, $f$ has finite order of growth.
We refer to \cite{Hay64,Jan85,Nev29,Nev53} for the terminology of Nevanlinna theory,
and in particular for the definitions of characteristic and order.

For entire functions we have stronger results.
Recall that an entire function $f$ is said to be of exponential type if the 
maximum modulus $M(r,f)$ satisfies $\log M(r,f)=\OO(r)$  as $r\to\infty$. 
The following result is a special case of a result  of
J.~Clunie and 
W.~K.~Hayman~\cite[Theorem~3]{Clu66}; see also~\cite[Theorem~4]{Pom}.

\begin{thm} \label{cluhay}
{\bf Clunie-Hayman Theorem.}
An entire function with bounded spherical derivative is of exponential type.
\end{thm}
\subsection{Further examples}\label{examples}
We have seen already some examples of properties for which Bloch's principle
holds, and we will see many more examples of such properties in the following
sections. But of course it is also very enlightening to consider counterexamples 
to the heuristic Bloch principle. The following examples are due to L.~A.~Rubel~\cite{Rub}.
\begin{exam}
{\bf Example.} 
Define $\langle f,D\rangle \in P$ if $f=g''$ for some function $g$ which
is holomorphic and univalent in $D$. Since the only univalent entire functions $g$ are
those of the form $g(z)=\alpha z+\beta$ where $\alpha ,\beta\in\C$, $\alpha \neq 0$, we see that 
$\langle f,\C\rangle \in P$ implies that $f(z)\equiv 0$. Thus $(a)$ holds.

On the other hand, 
define 
$$g_n(z):=n\left(z+\frac{1}{10}z^2+\frac{1}{10}z^3\right)$$ 
and $f_n:=g_n''$ for $n\in\N$.
Then 
$$\re g_n'(z)=n\left(1+\frac{1}{5}\re z +\frac{3}{10}\re z^2\right)>\frac{1}{2}>0$$
for $z\in \D$, which 
implies that $g_n$ is univalent in $\D$. So $\langle f_n,\D\rangle \in P$.
Now 
$f_n(z)=n(\frac15+\frac35 z)$ so that $f_n(-\frac13)=0$ 
while $f_n(0)\to\infty$ as $n\to\infty$.
Thus the $f_n$ do not form a normal family on $\D$. Hence $(b)$ fails.

It is easily seen that $P$ satisfies conditions $(i)$ and $(ii)$ of Zalcman's 
Principle, but condition $(iii)$ is not satisfied.
\end{exam}
\begin{exam}
{\bf Example.}
Define $\langle f,D\rangle \in P$ if $f$ is holomorphic in $D$ and 
$f'(z)\neq -1$, $f'(z)\neq -2$ and $f'(z)\neq f(z)$ for all $z\in D$.
If $\langle f,\C\rangle \in P$, then $f'$ is constant by Picard's Theorem so 
that  $f$ is   of the form $f(z)=\alpha z+\beta$ where $\alpha ,\beta\in\C$.
Since $f'(z)\neq f(z)$ we have $\alpha z+\beta-\alpha \neq 0$. This implies that $\alpha =0$ so that
$f$ is constant. Thus $(a)$ holds.

On the other hand, with $f_n(z):=nz$ we find that $\langle f_n,\D\rangle \in P$ 
for all $n\in\N$, but the $f_n$ do not form a normal family on $\D$ 
and thus $(b)$ fails.

Again it is clear that $(i)$ holds, and one can verify that $(iii)$ is also satisfied.
So in this case it is property $(ii)$ that fails to hold in Zalcman's 
Principle.
\end{exam}
\subsection{Other applications of Zalcman's Lemma}\label{other}
Not only can
Zalcman's lemma be used to prove for many properties $P$ that
statements $(a)$ and $(b)$ as above are {\em equivalent}, but it also often yields easy proofs 
that these statements are {\em true}. 

{\em We sketch a
simple proof of the theorems of Montel
and Picard}, due to A. Ros~\cite[p.~218]{Zal98}.
So let $\FF$ be as in the statement of Montel's Theorem and suppose that 
$\FF$ is not normal. Without loss of generality we may assume that 
$\{a_1,a_2,a_3\}=\{0,1,\infty\}$ and that $D$ is a disk. If $f\in \FF$
and $n\in\N$ there exists a function $g$ holomorphic in $D$ such that
$g^{2^n}=f$. Let $\FF_n$ be the family of all such functions $g$.
Note that
$$g^\#=\frac{1}{2^n}\frac{|f|^{1/2^{n}-1}|f'|}{1+|f|^{2/2^n}}
=\frac{1}{2^n}\frac{|f|^{-1}+|f|}{|f|^{-1/2^n}+|f|^{1/2^n}}f^\#
\geq \frac{1}{2^n}f^\#,$$
where we have used the inequality $a^{-1}+a\geq a^{-t}+a^t$ valid for 
$a>0$ and $0<t<1$. By Marty's Criterion, the family 
$\{f^\#:f\in \FF\}$ is not locally bounded. We deduce that, for 
fixed $n\in\N$, the family 
$\{g^\#:g\in \FF_n\}$ is not locally bounded. Using Marty's 
Criterion again we find that $\FF_n$ is not normal, for all $n\in\N$.

Note that if $g\in\FF_n$, then $g$ omits the values  $e^{2\pi i k/2^n}$
for $k\in {\mathbb Z}$. 
From the Zalcman Principle (and the remarks following it) we thus
deduce that there exists an entire function $g_n$  omitting
the values  $e^{2\pi i k/2^n}$ and satifying $g_n^\#(z)\leq g_n^\#(0)=1$.
The $g_n$ thus form a normal family and we have $g_{n_j}\to G$ for some 
subsequence $(g_{n_j})$ of $(g_n)$ and some nonconstant
entire function $G$. By Hurwitz's Theorem, $G$ omits the values 
$e^{2\pi i k/2^n}$ for all $k,n\in\N$. Since $G(\C)$ is open this 
implies that $|G(z)|\neq 1$ for all $z\in\C$. Thus either 
$|G(z)|< 1$ for all $z\in\C$ or $|G(z)|> 1$ for all $z\in\C$.  
In the first case $G$ is bounded and thus constant by Liouville's Theorem.
In the second case $1/G$ is bounded. Again $1/G$ and thus $G$ is constant.
Thus we get a contradiction in both cases.\qed 

The following result generalizes the theorems of Picard and Montel.
We have named it after R.~Nevanlinna (see~\cite[p.\ 102]{Nev29} or
\cite[\S X.3]{Nev53}) although only the
part concerning functions in the plane is due to him.  
The part concerning normal families is due to 
A.~Bloch~\cite[Theorem XLIV]{Blo26} and 
G.~Valiron~\cite[Theorem XXVI]{Val29}, with proofs
being based on Nevanlinna's theory, however.
A proof using different ideas was given by R.\ M.\ Robinson \cite{Rob39}.
For a proof using Zalcman's Lemma we refer to~\cite[\S 5.1]{Ber98}.
\begin{thm}\label{thmnevan}
{\bf Nevanlinna's Theorem.}
Let $q\in\N$, let 
$a_1,\dots,a_q\in\CC$ be distinct and let $m_1,\dots,m_q\in\N$.
Suppose that 
\begin{equation}\label{m_j}
\sum_{j=1}^q \left(1-\frac{1}{m_j}\right)> 2.
\end{equation}
Then the property that all $a_j$-points of $f$
have multiplicity at least $m_j$ is a Picard-Montel property.
\end{thm}
We may also allow $m_j=\infty$ here, meaning that $f$ has no $a_j$-points at all.

It is easy to see that the property $P$ in this result satisfies conditions $(i)$, $(ii)$ and
$(iii)$ of the Zalcman Principle~\ref{zalcprinc}
so that the statement about functions in the plane
is equivalent to that about normal families. But the proof that both statements are true
is more awkward.

{\em We sketch the proof of a special case of Nevanlinna's Theorem}. We
assume in addition that  the multiplicity $m(z)$ of each $a_j$-point $z$ of $f$ satisfies 
not only $m(z)\geq m_j$, but that in fact $m(z)$ is a multiple of $m_j$; that is,
$m(z)=n(z)m_j$ where $n(z)\in \N$. Let $P'$ be this property. Again $(i)$, $(ii)$ and
$(iii)$ are satisfied so that $P'$ is also a Bloch property. Moreover, properties $(a)$ and 
$(b)$ occuring in the definition of Bloch property are equivalent to property $(c)$
mentioned at the end of~\S\ref{seczhp}. Thus it suffices to prove $(c)$.

So let $\langle f,\C\rangle \in P'$ where $f$ has bounded spherical derivative, but suppose
that $f$ is not constant.
We may assume that $a_j\neq\infty$ for all $j$ and define
\begin{equation}\label{g(z)}
g(z)=\frac{f'(z)^M}{\prod_{j=1}^q(f(z)-a_j)^{(m_j-1)M/m_j}},
\end{equation}
where $M$ is the least common multiple of the $m_j$.
The assumption on the multiplicities of the $a_j$-points implies that $g$ does not 
have poles; that is, $g$ is entire.
Since $f$ is not constant there exists a sequence $(u_n)$ tending to $\infty$ 
such that $f(u_n)\to\infty$. The denominator of $g$ is a polynomial 
in $f$ of degree 
$$\ell:=\sum_{j=1}^q \frac{(m_j-1)M}{m_j}=M\sum_{j=1}^q 
\left(1-\frac{1}{m_j}\right)> 2M.$$
Thus 
$$\prod_{j=1}^q(f(u_n)-a_j)^{(m_j-1)M/m_j}\geq (1-o(1))|f(u_n)|^\ell
\geq (1-o(1))|f(u_n)|^{2M+1}.$$
On the other hand, 
$|f'(u_n)|\leq \OO(|f(u_n)|^2)$ as $n\to\infty$
since $f$ has bounded spherical derivative. 
Thus 
$|f'(u_n)|^M\leq o(|f(u_n)|^{2M+1})$ as $n\to\infty$. Overall we see that
$g(u_n)\to 0$. On the other hand, $g(z)\not\equiv 0$ since $f$ is not constant.
Thus $g$ is not constant, which implies that there exists
a sequence $(v_n)$ such that
$g(v_n)\to\infty$. It follows  that $f(v_n)\not\to\infty$.

We consider  $h_n(z):=f(z+v_n)$. Since $f$ has bounded spherical derivative, the
$h_n$ form a normal family.  Passing to a subsequence if necessary, we
may thus assume that $h_n$  converges locally 
uniformly to some meromorphic function $h:\C\to\CC$.
It follows that $h(z)\equiv a_k$ for some $k\in\{1,\dots,5\}$, because otherwise 
$$g(z+v_n)\to \frac{h'(z)^M}{\prod_{j=1}^q(h(z)-a_j)^{(m_j-1)M/m_j}}\neq\infty ,$$
which contradicts $g(v_n)\to\infty$.

Thus $h_n\to a_k$ as $n\to\infty$.
For sufficiently large $n$ the function 
$$\psi_n(z):=h_n(z)-a_k=f(z+v_n)-a_k$$
is holomorphic in $\D$,
and we have $\psi_n\to 0$ as $n\to\infty$.
Since the multiplicity of all $a_k$-points of $f$ is a multiple 
of $m_k$ we may define a holomorphic
branch $\phi_n$ of the $m_k$-th root of $\psi_n$; that 
is, $\phi_n:D\to \C$ and $\phi_n(z)^{m_k}=\psi_n(z)$.
We also have $\phi_n\to 0$. Thus $\phi_n'(z)\to 0$.
Now $(m_k|\phi_n'(z)|)^{m_k}=|\psi_n(z)|^{1-m_k}|\psi_n'(z)|^{m_k}$.
Hence $|f_n'(v_n)|^{m_k}/|f_n(v_n)-a_k|^{m_k-1}=
|\psi_n'(0)|^{m_k}/|\psi_n(0)|^{m_k-1}\to 0$.
This implies that $g(v_n)\to 0$, a contradiction.\qed

We mention that for $q=3$ the special case of Nevanlinna's Theorem proved above 
goes back to C.~Carath\'eodory~\cite{Car05}. Note that this special case still gives
a generalization of the theorems of Picard and Montel.

The general version of Nevanlinna's Theorem may be proved along the
same lines. In this case the functions $\phi_n$ occuring in the above proof may be
multi-valued. A similar argument as above may still me made, however, by using
a version of  Schwarz's Lemma for multivalued functions due to 
Z.~Nehari~\cite{Neh47};
see~\cite{Ber98} for more details.

We note that the hypothesis~(\ref{m_j}) is best possible. In fact,
if we have equality in~(\ref{m_j}), then there exists an elliptic function 
$f$ 
which  satisfies the differential equation
$f'(z)^M=\prod_{j=1}^q(f(z)-a_j)^{(m_j-1)M/m_j}$. 
It follows that 
$\langle f,\C\rangle\in P$ for $a$ nonconstant function $f$.  
Note that the 
function $g$ defined in the above proof by~(\ref{g(z)}) now satisfies $g(z)\equiv 1$.

The possible choices for the $m_j$ are $q=4$ and $m_j=2$ for all $j$ or
$q=3$ and $(m_1,m_2,m_3)=(3,3,3)$, $(m_1,m_2,m_3)=(2,3,6)$ or $(m_1,m_2,m_3)=(2,4,4)$,
up to permutations of the $m_j$.
One can modify the construction to allow the case that $a_j=\infty$ for some $j$.
Moreover, one can  also modify the above considerations to include the case
$m_j=\infty$. In this case the resulting functions $f$ are 
trigonometric functions or the exponential function, or obtained from these
functions by linear transformations.

Closely related to Nevanlinna's Theorem is one of the main results from the Ahlfors theory of
covering surfaces; see~\cite{Ahl35}, \cite[Chapter 5]{Hay64} or~\cite[Chapter XIII]{Nev53}.
Let $D\subset \CC$ be a domain and let
$f:D\to\CC$ be a meromorphic function. Let 
$V\subset\CC$ be a Jordan domain.
A simply-connected component $U$ of $f^{-1}(V)$ 
with $\overline{U}\subset D$ is called an {\em island}
of $f$ over $V$. 
Note that then $f|_U:U\to V$ is a proper map.
The degree of this proper map is called the 
{\em multiplicity} of the island $U$.
An island of multiplicity one is called a 
{\em simple island}.
\begin{thm}\label{ahl}
{\bf Ahlfors's Theorem.}
Let $q\in\N$,
$D_1,\dots,D_q\subset\CC$ Jordan domains 
with pairwise disjoint closures and $m_1,\dots,m_q\in\N$
satisfying $(\ref{m_j})$.
Then the property that all islands of $f$ over $D_j$
have multiplicity at least $m_j$ is a Picard-Montel property.
\end{thm}
 
We indicate how Zalcman's Lemma can be used to deduce 
Ahlfors's Theorem from Nevanlinna's Theorem; see~\cite{Ber98} for a more detailed 
discussion.

Fix $a_1,\dots,a_q\in\C$.
First we show that there exists $\varepsilon >0$ such that the conclusion of 
Ahlfors's Theorem is true if $D_j=D(a_j,\varepsilon)$. Here and in the following
$D(a,\varepsilon):=\{z\in\C:|z-a|<\varepsilon\}$ denotes the disk of radius $\varepsilon$
around a point $a\in\C$. 
If such an $\varepsilon $ does not exist, then we can choose a sequence $(\varepsilon_n)$ 
tending to $0$ and find a sequence $(f_n)$ of nonconstant 
meromorphic functions on $\C$
which have no island of multiplicity less than $m_j$ over $D(a_j,\varepsilon )$.
By the arguments of \S\ref{seczhp} 
we may assume that 
the $f_n$ have bounded spherical derivative and in fact that
$f_n^\#(z)\leq f_n^\#(0)=1$ for all $z\in\C$. It follows that the $f_n$ form a normal family,
and thus we may assume without loss of generality that $f_n\to f$ for some
meromorphic function $f$. Since $f^\#(0)=1$ we see that $f$ is not constant.
We also find that $f$ has no island of multiplicity less than $m_j$ over $D(a_j,\varepsilon )$,
for any $\varepsilon >0$.
But this implies that all $a_j$-points of $f$ have multiplicity at least $m_j$,
contradicting Nevanlinna's Theorem. Thus there exists $\varepsilon >0$ such 
that the conclusion of 
Ahlfors's Theorem holds if  $D_j=D(a_j,\varepsilon)$. 

In the second step we reduce the general case to this special case.
To do this we use quasiconformal mappings; see~\cite{LV} for an introduction
to this subject.
So suppose that $f$ is nonconstant and meromorphic in the plane
such that every island over $D_j$ has multiplicity at least $m_j$,
for $j\in\{1,\dots,q\}$. 
We note that there exists a quasiconformal map $\phi:\C\to\C$ with
$\phi(D_j)\subset D(a_j,\varepsilon)$ for $j\in\{1,\dots,q\}$,
and the quasiregular map $\phi\circ f$ can be factored as
$\phi\circ f=g\circ \psi$ with a nonconstant meromorphic function
$g:\C\to\CC$ and a quasiconformal map $\psi:\C\to\C$.
It then follows that every island of $g$ over $D(a_j,\varepsilon)$ 
has multiplicity at least $m_j$, contradicting the first step above.\qed

As mentioned, the theorems of Picard and Montel are a special case
of Nevanlinna's Theorem, namely the case $q=3$ and $m_1=m_2=m_3=\infty$. 
Applying the arguments used in the proof above to the theorems of Picard and Montel
instead of Nevanlinna's Theorem, we obtain a direct proof of the
following special case of Ahlfors's Theorem.
\begin{thm}\label{ahlspec}
{\bf Special case of Ahlfors's Theorem.}
Let
$D_1,\dots,D_3\subset\CC$ be Jordan domains 
with pairwise disjoint closures. 
Then the property that $f$ has no islands over any of the domains $D_j$
is a Picard-Montel property.
\end{thm}
Another important special case of Ahlfors's Theorem~\ref{ahl} is the
case $q=5$ and $m_j=2$ for all $j$. This case is known as the 
{\em Ahlfors Five Islands Theorem}.

\subsection{A modified Bloch principle} \label{modified}
Closely related to Montel's Theorem~\ref{thmon} and Picard's Theorem~\ref{thpic} is the
following result.
\begin{thm}\label{greatpic}
{\bf Great Picard Theorem.}
Let $a_1,a_2,a_3\in \CC$ be distinct,  $D\subset \C$ a domain, $\zeta\in D$ and 
$f:D\setminus\{\zeta\}\to\CC$ meromorphic.
If $f(z)\neq a_j$ for all $j\in\{1,2,3\}$ and 
all $z\in D\setminus\{\zeta\}$, then $\zeta$ is not an essential singularity of $f$.
\end{thm}
This results suggests
a modification of the Bloch principle which says that for 
a Picard-Montel property $P$  there should not exist a meromorphic function
having the property $P$ in the neighborhood of an essential 
singularity. 
In other words, for a property $P$
of meromorphic functions the conditions $(a)$, $(b)$ and $(c)$ discussed
in~\S\ref{seczhp} should imply the following condition: 
\begin{itemize}
\item[$(d)$] {\em If $\langle f,\D\setminus\{\zeta\}\rangle \in P$ for some domain $D$ 
and some $\zeta\in D$, then $\zeta$ is not essential singularity of $f$.}
\end{itemize}
D.~Minda~\cite{Min85} gives a discussion of this modification of the heuristic
principle, and he shows that for holomorphic families this 
modified heuristic principle
holds under the hypotheses of Zalcman's Principle~\ref{zalcprinc}.
\begin{thm}\label{thmminda}
{\bf Minda's Principle.}
 Suppose that a property $P$ of holomorphic functions satisfies 
the conditions $(i)$, $(ii)$ and $(iii)$ of 
Zalcman's Principle~\ref{zalcprinc}.
Then each of the conditions $(a)$, $(b)$ and  $(c)$ introduced
in~\S\ref{seczhp} implies that $(d)$ holds.
\end{thm}
Minda also points out that there are 
{\em meromorphic} families
where Zalcman's Principle applies, but where the modified heuristic
principle does not hold.
However, if one adds a further condition to 
the ones given by Zalcman, then the modified heuristic principle
holds; see~Minda~\cite[\S 5]{Min85} and~\cite[\S 3]{Ber00} for the following result.
\begin{thm} \label{modexp}
{\bf Theorem.} Suppose that a property $P$ of meromorphic functions satisfies 
the conditions $(i)$, $(ii)$ and $(iii)$ of Zalcman's Principle~\ref{zalcprinc}.
Suppose that $P$ satifies in addition the condition 
\begin{itemize}
\item[$(iv)$] If $\langle f,\C\setminus\{0\}\rangle \in P$, then
$\langle f\circ \exp,\C\rangle \in P$. 
\end{itemize}
Then each of the conditions $(a)$, $(b)$ and  $(c)$ implies that $(d)$ holds.
\end{thm}
The additional condition $(iv)$ can be compared with $(ii)$. Both conditions are 
obviously satisfied for properties $P$ which concern only the range of $f$. 
It is easily seen that the properties occuring in Nevanlinna's and Ahlfors's Theorem
in~\S\ref{other} satisfy~$(iv)$.

The proof of 
these theorems
uses the following results due to 
O.~Lehto and
L.~V.\ Virtanen~\cite{Leh57,Leh59,LV57}.
\begin{thm} \label{lehto}
{\bf Lehto-Virtanen Theorem.}
Suppose that a meromorphic function $f$ has
an essential singularity at $\zeta$. Then
\[
\limsup_{z\to\zeta}|z-\zeta| f^\#(z) \geq\frac{1}{2}.
\]
If $f$ is holomorphic, then
\[
\limsup_{z\to\zeta}|z-\zeta| f^\#(z) =\infty.
\]
\end{thm}
Lehto and Virtanen~\cite{LV57}
had shown that 
$\limsup_{z\to\zeta}|z-\zeta| f^\#(z)\geq k$ for some absolute constant
$k>0$. This weaker result would suffice for our purposes.
Lehto~\cite{Leh57} later showed that one can take $k=\frac12$.

{\em Proof of Minda's Principle~\ref{thmminda}}.
Since we assume that the hypotheses of 
Zalc\-man's Principle~\ref{zalcprinc} are satisfied,
the conclusion of~\ref{zalcprinc} also holds. Thus 
conditions
$(a)$ and  $(b)$ are equivalent, and the discussion after~\ref{zalcprinc} 
shows that these conditions are also equivalent to $(c)$.
Suppose now that one and hence 
all of these conditions are satisfied. We want to prove $(d)$.

So let
$\langle f,D\setminus\{\zeta\}\rangle \in P$, where $D$ is a domain, $f$ is holomorphic
and $\zeta\in D$.  Suppose that $\zeta$ is an essential singularity of $f$.
We may assume that $\zeta=0$. 
By the Lehto-Virtanen Theorem
there exists a sequence $(c_n)$ in $D$
such that $c_n\to 0$ and
$|c_n|f^\#(c_n)\to\infty$. 
For sufficiently large $n$ the function
$f_n(z):=f(c_n + c_n z)$ is then holomorphic in the unit disk $\D$,
and $f_n^\#(0) =c_n f^\#(c_n)\to\infty$.
Thus the $f_n$ cannot form a normal family by Marty's Criterion.
On the other hand, we deduce from~$(i)$ and~$(ii)$ that the $f_n$ also satisfy $P$.
This is a contradiction to~$(b)$.\qed

{\em To prove Theorem~\ref{modexp}}, we note again that
$(a)$, $(b)$ and $(c)$ are equivalent.  
We suppose that these conditions are satisfied and  want to prove $(d)$.

So let $\langle f,D\setminus\{\zeta\}\rangle \in P$ for a domain $D$ 
and some $\zeta\in D$. We may assume that $\zeta=0$. 
Suppose that $\zeta$ is an essential singularity.
By the Lehto-Virtanen Theorem
there exists sequence $(c_n)$ in $D$
such that $c_n\to 0$ and
$|c_n|f^\#(c_n)\geq \frac14$. 
Choose $r>0$ such that $D(0,r)\subset D$.
We define $r_n:=r/|c_n|$ 
and 
$g_n:D(0,r_n)\to \CC$ by
$g_n(z):=f(c_n z)$. 
Since $P$ satisfies condition $(b)$, as well as $(i)$ and $(ii)$,
and since $r_n\to\infty$, we see that the $g_n$ form a normal family.
Without loss of generality we may assume that $g_n\to g$ for some 
$g:\C\backslash \{0\}\to \CC$. Since
$g_n^\#(1)=c_nf^\#(c_n)\geq \frac14$ we have $g^\#(1)\geq \frac14$
so that $g$ is nonconstant. It follows from $(iii)$ that 
$\langle g,\C\setminus\{0\}\rangle \in P$. By condition $(iv)$ thus 
$\langle g\circ \exp,\C\rangle \in P$. From $(a)$ we 
may deduce that $g\circ \exp$
is constant, and so is $g$, a contradiction.\qed 
\subsection{Quasinormality}
We cannot expect that condition $(d)$ of \S\ref{modified}
implies the conditions $(a)$, $(b)$, and $(c)$ introduced in~\S\ref{seczhp}.
For example, the condition that $f$ takes three values $a_1,a_2,a_3$ only 
$N$ times, for some fixed  number $N\in\N$, satisfies condition $(d)$, but
none of the conditions $(a)$, $(b)$, and $(c)$.  This condition does, however,
imply {\em quasinormality}.

We say that 
a family $\FF$ of functions meromorphic in a domain
$D$ is {\em quasinormal} 
(cf.\ \cite{Chu,Mon27,Schi})
if for each sequence $(f_k)$ in $\FF$ there exists a 
subsequence $(f_{k_j})$ and a finite set $E\subset D$ such that
$\left(f_{k_j}\right)$ converges locally uniformly in
$D\backslash E$.
If the cardinality of the exceptional set $E$ can be bounded
independently of the sequence $(f_k)$, and
if $q$ is the smallest 
such bound,
then we
say that $\FF$ is quasinormal of {\em order}~$q$.

Many of the results about normal families have extensions involving the concept of
quasinormality. Here we only mention the corresponding generalization of 
Montel's Theorem, also proved by Montel~\cite[p.~149]{Mon27}.
\begin{thm}
{\bf Montel's Theorem.}
Let $a_1,a_2,a_3\in \CC$ be distinct, 
let $0\leq m_1\leq m_2\leq m_3$, 
let $D\subset \C$ be a domain
and let $\FF$ be a family of functions meromorphic
in $D$.
Suppose  that $f$ takes the value $a_j$ at most $m_j$ times in $D$,
for all $j\in\{1,2,3\}$ and all $f\in\FF$. Then $\FF$ is quasinormal
of order at most $m_2$.
\end{thm}
Quasinormality will also be discussed in~\S\ref{quasinorm} below.
A detailed study of quasinormality, and in fact of a more general 
concept called $Q_m$-normality, has been given by C.~T.~Chuang~\cite{Chu}.

\section{Exceptional values of derivatives}\label{exceptional}
\subsection{Introduction}\label{intro2}
We discuss some variants of the theorems of Picard and Montel where
exceptional values of $f$ are replaced by those of a derivative.
Our starting point is the following result proved by
G.\ P\'olya's student W.\ Saxer in~1923; see~\cite[Hilfssatz, p.~210]{Sax}
and~\cite{Sax26}.
\begin{thm}\label{saxerthm}
{\bf Saxer's Theorem.}
Let $f$ be a transcendental entire function and let $a,b\in\C$. Suppose that
the equations $f(z)=a$ and $f'(z)=b$ have only finitely many solutions.
Then $b=0$.
\end{thm}
Combining this with the Great Picard Theorem, applied 
to $f'$, we see that if a transcendental entire function
$f$ takes a value $a$ only finitely many times, then $f'$ takes every
nonzero value infinitely often. In 1929, E.\ Ullrich~\cite[p.~599]{Ull}
showed that under this hypothesis all derivatives $f^{(k)}$, $k\geq 1$,
take every nonzero value infinitely often.
Ullrich writes that this result has been known for several years,
and he attributes it to P\'olya and Saxer~\cite{Sax},
although it can be found there only for the first derivative.

A simple discussion of the case where $f$ is a polynomial
now leads to the following result.
\begin{thm}\label{psu}
{\bf P\'olya-Saxer-Ullrich Theorem.}
Let $f$ be an entire function and let $k\geq 1$. Suppose that
$f(z)\neq 0$ and $f^{(k)}(z)\neq 1$ for all $z\in\C$. Then $f$ is constant.
\end{thm}
Here and in the following theorems the conditions
$f(z)\neq 0$ and $f^{(k)}(z)\neq 1$ can  be replaced by 
$f(z)\neq a$ and $f^{(k)}(z)\neq b$ as long as $b\neq 0$.

Ullrich's result was obtained independently a few years later by F.\
Bureau~\cite{Bur31,Bur32,Bur97}. 
Actually Bureau considered functions with an essential singularity,
as in~\S\ref{modified}. Bureau also gave a normality result, but
he required additional conditions 
besides $f(z)\neq 0$ and $f^{(k)}(z)\neq 1$.
The complete normal family analogue was obtained by
C.\ Miranda~\cite{Mir} in~1935.
\begin{thm} \label{thmburmir}
{\bf Miranda's Theorem.}
Let $\FF$ be a family of  functions holomorphic in a domain $D$ and $k\geq 1$. Suppose that
$f(z)\neq 0$ and $f^{(k)}(z)\neq 1$
for all $f\in\FF$ and $z\in D$. Then $\FF$ is normal.
\end{thm}

In 1959, W.~K.~Hayman~\cite[Theorem~1]{Hay59} extended the 
P\'olya-Saxer-Ullrich Theorem to meromorphic functions.
We mention 
that the case that $f$ is meromorphic had also been considered by Ullrich,
but in this case he required additional hypotheses, e.~g. that~$\infty$
is a Borel exceptional value~\cite[p.~599]{Ull}.
\begin{thm}\label{thhayman}
{\bf Hayman's Theorem.}
Let $f$ be meromorphic in the plane and $k\geq 1$. Suppose that
$f(z)\neq 0$ and $f^{(k)}(z)\neq 1$
for all $z\in\C$. Then $f$ is constant.
\end{thm}
\begin{rem}
More generally, Hayman proved that if $f$ and $f^{(k)}-1$ have only finitely many zeros,
then $f$ is rational.
\end{rem}
It took 20 years until Y.\ Gu~\cite{Gu} proved the normal family analogue of Hayman's Theorem.
\begin{thm}
{\bf Gu's Theorem.}
Let $\FF$ be a family of  functions meromorphic in a domain $D$ and $k\geq 1$. Suppose that
$f(z)\neq 0$ and $f^{(k)}(z)\neq 1$
for all $f\in\FF$ and $z\in D$. Then $\FF$ is normal.
\end{thm}

\subsection{A generalization of Zalcman's Lemma}\label{genzal}
Zalcman's Principle as stated in~\S\ref{seczhp} does not 
apply to conditions such as 
``$f(z)\neq 0$ and $f^{(k)}(z)\neq 1$.'' 
However, 
it was shown by X.~Pang~\cite{Pang89,Pang90}
that there is an extension of Zalcman's Lemma which allows 
to deal with such conditions.
We state this extension in its most general form, and not only 
in the form needed to prove that
``$f(z)\neq 0$ and $f^{(k)}(z)\neq 1$'' is a Bloch property.
\begin{thm}
{\bf Zalcman-Pang Lemma.}
Let $\FF$ be a family of functions meromorphic in a domain 
$D\subset\C$ and let $m\in\N$, $K\geq 0$ and  $\alpha\in\R$ 
with $-m\leq\alpha< 1$.
Suppose that the zeros of the functions in $\FF$ have 
multiplicity at least $m$; that is, if $f\in\FF$ and $\xi\in D$
with $f(\xi)=0$, then $f^{(k)}(\xi)=0$ for $1\leq k\leq m-1$.
If $\alpha=-m$, then suppose in addition that 
$|f^{(m)}(\xi)|\leq K$ if $f\in\FF$, $\xi\in D$ and $f(\xi)=0$.

Suppose that $\FF$ is not normal at $z_0\in D$.
Then there exist a sequence
$(f_k)$ in $\FF$, a sequence $(z_k)$ in $D$, a sequence $(\rho_k)$
of positive real numbers and a nonconstant meromorphic function $f:\C\to\CC$ 
 such that $z_k\to z_0$,
$\rho_k\to 0$ and
$$\rho_k^\alpha f_k(z_k+\rho_k z)\to f(z)$$
locally uniformly in $\C$.
Moreover,  
$f^\#(z)\leq f^\#(0)=mK+1$ for all $z\in \C$.
\end{thm}
We omit the proof, but mention only that 
its spirit is similar to the 
proof 
of the original Zalcman Lemma, but the technical details are
more convoluted.

We note that $\{1/f:f\in\FF\}$ is normal if and only if
$\FF$ is normal.
Thus we obtain an analogous result for $-1<\alpha\leq\ell$, if
the poles of the functions in $\FF$ have multiplicity at least 
$\ell$, with an additional hypotheses if $\alpha=\ell$.
Note that no hypothesis on the zeros or poles is required when
$-1<\alpha<1$.

The Zalcman Lemma~\ref{lemzalc} is of course the case $\alpha=0$.
As mentioned,
the idea to introduce the exponent $\alpha$ seems to be due to 
X.~Pang~\cite{Pang89,Pang90}, who proved that one can always 
take $-1<\alpha<1$.
It was shown by X.~Pang and G.~Xue~\cite{PX} that 
$\alpha<0$ is admissible if the functions in $\FF$ 
have no zeros.
Then H.~Chen and X.~Gu~\cite[Theorem~2]{CG}
proved that one can take $-m<\alpha\leq 0$ if the zeros
of the functions in $\FF$ have multiplicity at least $m$.
Finally, the case $\alpha=-m$ is due to 
X.~Pang and L.~Zalc\-man~\cite[Lemma~2]{PZ}.
The special case $\alpha=-m=-1$ had been treated before 
by X.~Pang~\cite{Pang02}.

The Zalcman-Pang Lemma shows that the Zalcman Principle~\ref{zalcprinc}
may be modified by 
replacing the condition $(ii)$ by 
\begin{itemize}
\item[$(ii')$] {\em There exists $\alpha\in (-1,1)$ such that 
if $\langle f,D\rangle \in P$ and $\varphi(z)=\rho z +c$ where $\rho,c\in\C$, $\rho\neq 0$,
then $\langle \rho^\alpha (f\circ \varphi),\varphi^{-1}(D)\rangle \in P$.}
\end{itemize}
or 
\begin{itemize}
\item[$(ii'')$] {\em There exists $m\in\N$ and $\alpha\in (-m,1)$ such that 
if $\langle f,D\rangle \in P$, then all zeros of $f$ have multiplicity at least $m$,
 and if $\varphi(z)=\rho z +c$ where $\rho,c\in\C$, $\rho\neq 0$,
then $\langle \rho^\alpha (f\circ \varphi),\varphi^{-1}(D)\rangle \in P$.}
\end{itemize}
We leave it to the reader to formulate a condition for the case that $\alpha=-m$,
or the case that the functions with property $P$ have only multiple poles of
multiplicity at least $\ell$.

If $\langle f,D\rangle \in P$ implies that $f$ has no zeros, then we can take
any $m$ and $\alpha$ in~$(ii'')$. The choice $\alpha=-k$ implies that the
property occuring in the theorems of 
P\'olya-Saxer-Ullrich, Miranda, Hayman and Gu is indeed
a Bloch property. In particular, Gu's Theorem  can then 
be deduced from Hayman's, and Miranda's from that of P\'olya-Saxer-Ullrich.
However, as with the theorems of Picard and Montel, the Zalcman-Pang Lemma not only
shows that the theorems of Miranda and P\'olya-Saxer-Ullrich
are {\em equivalent}, it can 
also be used to {\em prove} them; see~\cite{Grahl,PZnewzealand}.

{\em Proof
of Miranda's theorem (following~\cite{PZnewzealand}) }.
We write $\langle f,D\rangle\in P$ if $f$ is holomorphic in $D$ and if 
$f$ and $f^{(k)}-1$ have no zeros in $D$. As seen above, it follows from the
Zalcman Principle, with $(ii)$ replaced by $(ii'')$, that $P$ is 
a Bloch property. The argument given at the end of~\S\ref{seczhp} shows that 
the theorems of P\'olya-Saxer-Ullrich and Miranda are equivalent to the condition 
$(c)$ given there.

Let thus $\langle f,\C\rangle\in P$ where $f$ has bounded spherical derivative.
By the Clunie-Hayman Theorem~\ref{cluhay},
$f$ has exponential type. Since $f$ has no zeros 
this implies that $f$ has the form $f(z)=e^{az+b}$ with $a,b\in\C$.
Since $f^{(k)}(z)-1=a^ke^{az+b}-1$ has no zeros we deduce that $a=0$ and
thus $f$ is constant.\qed

We mention, however, that we have not been able to find a proof of Hayman's or
Gu's Theorem based on these ideas.

\subsection{Multiple values of derivatives}\label{multiple}
Nevanlinna's Theorem~\ref{thmnevan}
can be seen as a generalization of the 
theorems of Picard and Montel, where the hypothesis that $f$ does not take 
a value $a_j$ is replaced by the hypothesis that the $a_j$-points of $f$ have
high multiplicity. One may ask whether the theorems of Hayman and Gu (or of 
P\'olya-Saxer-Ullrich
and Miranda) admit similar generalizations. This question was considered 
in~\cite{BL1,Chen,Dra69,WangFang,YZ1,YZ2}. The results below are taken mostly
from~\cite{BL1}.

Let $k$ be a positive integer and let
$0 < M \leq \infty $,
$0 < N \leq \infty $. For a function $f$ meromorphic in a domain $D$ we
say that $f$ has the property $P(k,M,N)$, written again as
$\langle f,D\rangle \in P(k,M,N)$, if all zeros of $f$ in $D$ have 
multiplicity at least $M$, while
all zeros of $f^{(k)}-1$ in $D$ have multiplicity at least $N$.
Here  $M = \infty $ or $ N = \infty $
should be interpreted as meaning that there
are no corresponding zeros in $D$. 
\begin{thm}\label{thm0}
{\bf Theorem.} 
Let $k, M, N \in \N$. Then $P(k,M,N)$ is a Bloch property.
\end{thm}
{\em Proof}. Suppose first that $M\leq k$.
Define
$f_n (z) = 2n(z - a)^k$ with  
$n \in \N$ and 
$a \in \C$. 
Then $f_n$ satisfies $P(k,M,N)$ for any $N$. Moreover,
$f_n$ is nonconstant and entire, and there
is no neighbourhood of $a$ on which the
$f_n$ form a normal family. So both statements $(a)$ and $(b)$ occuring in the
definition of Bloch property in~\S\ref{seczhp} are false.

Suppose now that $M>k$. Then $P(k,M,N)$ satisfies condition $(ii'')$ in~\S\ref{genzal}
with $m:=M$ and $\alpha:=-k$. It is obvious that $P(k,M,N)$ also satisfies
condition~$(i)$ of Zalcman's
Principle~\ref{zalcprinc}. An application of 
Hurwitz's Theorem shows that condition~$(iii)$ of~\ref{zalcprinc} is
also satisfied. The conclusion thus follows from the Zalcman Principle,
as generalized in~\S\ref{genzal}.
\qed

Theorem~\ref{thm0} says that whether a family 
$\{f: \langle f,D\rangle \in P(k,M,N)\}$ is normal is equivalent to whether 
$\langle f,\C\rangle \in P(k,M,N)$ contains only constant functions~$f$.
It does not say whether these statements are true or false.
For example, for $(k,M,N)=(1,3,3)$ both statements are false, as shown by the example 
$f:=1/\wp'$, where $\wp$ is the Weierstra\ss\ 
elliptic
function satisfying 
the equation
$(\wp')^2 = 4 \wp^3 - g_2 \wp  - g_3 
=  \wp^3 + 3 \wp - 1$.
Then $f$ has only triple zeros and so has
$f'-1=-4(\wp+\frac12)3/(\wp')^2$. Thus $\langle f,\C\rangle \in P(1,3,3)$,
and $f$ is nonconstant.
But I do not know whether property $P(k,M,N)$ forces
a meromorphic function in the plane to be constant -- and a family of 
meromorphic functions
to be normal -- if, for example,  $(k,M,N)=(1,3,4)$ or $(k,M,N)=(1,4,3)$.

So while the precise conditions on $M$ and $N$ yielding this are not known, 
some partial answers are available.
\begin{thm}\label{thm1}
{\bf Theorem.} 
Let $k$ be a positive integer and let
$0 < M \leq \infty $ and
$0 < N \leq \infty $ with
$$
\frac{2k + 3 + 2/k}{M} +
\frac{2k +4+2/k }{N} < 1.
\label{h2}
$$
Then $P(k,M,N)$ is a Picard-Montel property.
\end{thm}
We omit the proof, which -- among other things -- is based on  Nevanlinna theory. The interested
reader is referred to~\cite{BL1}.

The following result is due to Y.~Wang and M.~Fang~\cite[Theorem~7]{WangFang}.
\begin{thm}
{\bf Wang-Fang Theorem.}
Let $k$ be a positive integer. 
Then $P(k,k+2,\infty)$ is a Picard-Montel property.
\end{thm}
One ingredient in the proof of this theorem
is the following result proved in~\cite[Corollary~3]{BE}.
Here a complex number $w$ is called a critical value of $f$ if there exists 
$\zeta$ such that $f'(\zeta)=0$ and $f(\zeta)=w$. 
\begin{thm}\label{critasy}
{\bf Theorem.} 
If a meromorphic function of finite order $\rho$ has only
finitely many critical values, then it has at most $2\rho$
asymptotic values.
\end{thm}
{\em We give a
 sketch of the proof of Wang-Fang Theorem}. First
note that -- as remarked at the end of~\S\ref{seczhp} --
it suffices to prove condition $(c)$ for the property $P=P(k,k+2,\infty)$.
So suppose that $f$ is a function meromorphic in the plane which has
bounded spherical derivative such that all zeros of $f$ have multiplicity $k+2$
at least and such that $f^{(k)}-1$ has no zeros. As remarked in~\S\ref{seczhp},
$f$ has finite order.

We consider the auxiliary function $g(z)=z-f^{(k-1)}(z)$. 
Then $g'$ has no zeros. Moreover, $g$ has finite order. 
By Theorem~\ref{critasy}, $g$ has only finitely many asymptotic values.

Suppose now that $f$ has a zero $\zeta$. Since this zero has multiplicity $k+2$ 
at least, we find that $g(\zeta)=\zeta$ and $g'(\zeta)=1$. 
In the terminology of complex dynamics
(see, for example,~\cite{Bea,Ber93,Mil99,Ste}) 
the point $\zeta$ is thus a parabolic fixed point of $g$.
By a 
classical result from
complex dynamics, sometimes called the Leau-Fatou Flower Theorem, 
there exists a domain $U$ 
with $\zeta\in\partial U$
where the iterates 
$g^n$ of $g$ tend  to~$\zeta$ as $n\to\infty$.
Moreover, a maximal domain $U$ 
with this property contains a critical or asymptotic value of $g$. 
Since $g$ has no critical and only finitely many asymptotic values,
this implies that $f$ has only finitely many zeros.
Hayman's Theorem~\ref{thhayman}, or more
precisely the remark made after it,
implies that $f$ is rational. A discussion of this case then completes the proof.\qed

The examples 
\begin{equation} \label{exrat}
f(z) = \frac{(z + a)^{k+1}}{k! (z+b)}
\end{equation}
with $a, b \in \C, a\neq b$ show that 
$P(k,k+1,\infty)$ is not a Picard-Montel property.
However, Wang and Fang~\cite[Lemma 8]{WangFang}
have shown that every nonconstant rational function $f$ satisfing 
$\langle f,\C\rangle \in P(k,k+1,\infty)$ has this form. 
The argument used in the proof of the Wang-Fang Theorem shows that 
there is no transcendental function $f$ meromorphic in the plane 
and of finite order such that 
$\langle f,\C\rangle \in P(k,k+1,\infty)$. 
However, this restriction on the order turns out not to be
necessary.
In fact, 
Wang and Fang~\cite[Theorem 3]{WangFang} proved that 
there is no transcendental  function $f$ meromorphic in the plane
satisfying $\langle f,\C\rangle \in P(k,3,\infty)$ for some $k\in\N$.
And Nevo, Pang and Zalcman~\cite{NPZ2} have recently shown that 
there is no transcendental function $f$ meromorphic in the plane 
such that $\langle f,\C\rangle \in P(1,2,\infty)$. 

The examples (\ref{exrat}) also show that the family 
$\{f: \langle f,D\rangle \in P(k,k+1,\infty)\}$
is not normal for any domain $D$ and $k\in\N$. However, 
Nevo, Pang and Zalcman~\cite{NPZ2} have recently shown that
$\{f: \langle f,D\rangle \in P(1,2,\infty)\}$
is quasinormal of order $1$,
and, as Larry Zalcman has kindly informed me, their
method can be extended to
yield that $\{f: \langle f,D\rangle \in P(k,k+1,\infty)\}$
is quasinormal of order $1$ for every $k\in\N$.
For further results in this direction we refer to~\cite{NP,NPZ}.

The following result strengthens the Wang-Fang Theorem.
\begin{thm}\label{thm5}
{\bf Theorem.}
Let $k$ be a positive integer. Then there exists a positive
integer $T_k$ such that 
$P(k,k+2,T_k)$ is a Picard-Montel property.
\end{thm}
{\em Proof}.
Again it suffices that to show that there exists $T_k$ such that $(c)$ holds
for the property $P(k,k+2,T_k)$. 
Suppose that this is not the case. 
Then for each $n\in\N$ there exists a nonconstant entire function $f_n$ such
that $\langle f_n,\C\rangle \in P(k,k+2,n)$, and we may assume that 
$f_n^\#(z)\leq f_n^\#(0)=1$ for all $z\in\C$. 
Thus the $f_n$ form a normal family so that 
$f_{n_j}\to f$ for some subsequence $(f_{n_j})$ of $(f_n)$.
Hurwitz's Theorem now implies that $\langle f,\C\rangle \in P(k,k+2,\infty)$, 
contradicting the Wang-Fang Theorem.\qed
\subsection{Exceptional values of differential polynomials}\label{diffpol}
Let $n\in\N$ and $a,b\in\C$, $a\neq 0$. 
We say that a meromorphic function $f:D\to\C$ satisfies
property $P(n,a)$ if $f(z)^nf'(z)\neq a$  for all $z\in D$. 
And we say that
$f$ has
property $Q(a,b,n)$ if $f'(z)+af(z)^n\neq b$ for all $z\in D$. 
Note that 
$\langle f,D\rangle \in P(n,a)$ if and only if 
$\langle 1/f,D\rangle \in Q(n+2,a,0)$.

It follows from the Zalcman-Pang Lemma with $\alpha:=-1/(n+1)$ that $P(n,a)$ is 
a Bloch property. This was a key ingredient in the proof of the following 
result~\cite{BE,CF,Zal94}.
\begin{thm}\label{thm6}
{\bf Theorem.}
$P(n,a)$ is a Picard-Montel property for all $n\geq 1$ and $a\neq 0$.
\end{thm}
That functions meromorphic in the plane which satisfy $P(n,a)$ are constant was
proved by Hayman
\cite[Corollary to Theorem 9]{Hay59} for  $n\geq 3$
and by E.\ Mues~\cite[Satz 3]{Mue79} for $n=2$.
The case that $f$ is entire is due to W.~K.\ Hayman 
\cite[Theorem 10]{Hay59} if $n\geq 2$ and to 
J.\ Clunie~\cite{Cl} if $n=1$.
By the remarks made above, Theorem~\ref{thm6} is equivalent to the following
result.
\begin{thm}\label{thm7}
{\bf Theorem.}
$Q(n,a,0)$ is a Picard-Montel property for all $n\geq 3$ and $a\neq 0$.
\end{thm}
Hayman
\cite[Theorem 9]{Hay59} also proved the following result.
\begin{thm}\label{thm8}
{\bf Theorem.}
If $\langle f,\C\rangle \in Q(n,a,b)$ where $a,b\in\C$, $a\neq 0$, and $n\geq 5$,
then $f$ is constant.
\end{thm}
The conclusion of this theorem is not true for $n=3$ and $n=4$, as shown by 
examples due to Mues~\cite{Mue79}.
For $n=4$ such an example is given
by $f(z):=\tan z$. Then $f'(z)=1+f(z)^2\neq 0$ so that
 $f'+\frac12 f^4-\frac12 =\frac12 (1+f^2)^2$ has no zeros.
Thus $\langle \tan,\C\rangle \in Q\left(4,\frac12,\frac12\right).$
The examples for $n=3$, or for different values of $a$ and $b$, are similar.

However, an argument due to Pang~\cite{Pang89,Pang90} shows that $Q(n,a,b)$ does
imply normality for $n\geq 3$. 
\begin{thm}\label{thmpang}
{\bf Pang's Theorem.}
If $n\geq 3$ and $a\neq 0$, then 
the family $\{f: \langle f,D\rangle \in Q(n,a,b)\}$ is normal on $D$ for 
each domain $D\subset \C$.
\end{thm}
{\em We briefly explain the idea of the proof}, which is by reduction to Theorem~\ref{thm7}.
We note that $Q(n,a,b)$ does not satisfy the condition $(ii')$ stated after
the Zalcman-Pang Lemma in~\S\ref{genzal}. However, 
with $\alpha:= 1/(n-1)\in (-1,1)$ and $\varphi(z):=\rho z +c$ where $\rho,c\in\C$, $\rho\neq 0$
we see that 
$\langle f,D\rangle \in Q(n,a,b)$ 
is equivalent to 
$$\langle \rho^\alpha (f\circ \varphi),\varphi^{-1}(D)\rangle \in Q(n,a,\rho^{n/(n-1)}b).$$
Since in the Zalcman-Pang Lemma one considers a sequence of $\rho$-values tending
to $0$ we see that the limit function $f$ occuring in this lemma satisfies
$\langle f,\C\rangle \in Q(n,a,0)$, contradicting Theorem~\ref{thm7}.\qed

When Pang wrote his papers, the conclusion of  Theorem~\ref{thm7} was known only 
for $n\geq  4$. Therefore he could proof his  result only for $n\geq 4$. But his
argument also extended to the case $n=3$, once
Theorem~\ref{thm7} was known; see also~\cite[p.~143]{Schi} for further discussion.

For $n\in\{3,4\}$ and $a,b\neq 0$ the property $Q(n,a,b)$ 
is thus a counterexample to the Bloch principle as stated in~\S\ref{seczhp}.
However, 
L.~Zalcman and, independently, the present writer
have in recent years suggested a variant of Bloch's principle 
which allows to deal with properties such as $Q(n,a,b)$.
Instead of the condition 
\begin{itemize}
\item[$(a)$] {\em If $\langle f,\C\rangle \in P$, then $f$ is constant.}
\end{itemize}
introduced in \S\ref{seczhp}
we consider the following condition:
\begin{itemize}
\item[$(a')$] {\em The family $\{f: \langle f,\C\rangle \in P\}$ is normal on $\C$.}
\end{itemize}
Note that $(a')$ is satisfied in particular if 
$\{f: \langle f,\C\rangle \in P\}$ 
consists only of constant functions.
In other words, the condition $(a)$ 
implies $(a')$. Recall the condition $(b)$ in the formulation 
of the original Bloch principle:
\begin{itemize}
\item[$(b)$] {\em The family $\{f: \langle f,D\rangle \in P\}$ is normal on $D$ for 
each domain $D\subset \C$.}
\end{itemize}
The variant of Bloch's principle mentioned
now says that $(a')$ should be equivalent to $(b)$.
Note that $(b)$ trivially implies $(a')$, so what this modification
of Bloch's principle
is really asking for is that $(a')$ implies $(b)$.

We have seen that for $n\in\{3,4\}$ and $a,b\neq 0$ the property $P:=Q(n,a,b)$ 
is an example where the original Bloch principle ``$(a)\Leftrightarrow (b)$''
fails, but where the variant ``$(a')\Leftrightarrow (b)$'' holds.

Some further cases where this is true will be discussed in
\S\ref{omitzero}, \S\ref{intro3} and \S\ref{sharedvalues}.

\subsection{Meromorphic functions with derivatives omitting zero}\label{omitzero}
The following result was proved by 
W.~K.\ Hayman~\cite[Theorem~5]{Hay59} for $k=2$ and
by
J.~Clunie~\cite{Clu62}
for $k\geq 3$.
\begin{thm}\label{thm11}
{\bf Hayman-Clunie Theorem.}
Let $f$ be entire
and let $k\geq 2$.
Suppose that $f$ and $f^{(k)}$ have no
zeros. Then $f$ has the form
$f(z)=e^{a z+b}$ 
where $a,b\in\C$, $a\neq 0$. 
\end{thm}
Hayman obtained the case $k=2$ as a corollary of his Theorem~\ref{thhayman}.
In fact, if $f$ and $f''$ have no zeros, then $F:=f/f'$ satisfies $F(z)\neq 0$
and $F'(z)-1=-f(z)f''(z)/f'(z)^2\neq 0$ for all $z\in\C$. Thus $F$ is constant
by Theorem~\ref{thhayman}, and this implies that $f$ has the form stated.

So we see that the conclusion of the above theorem is equivalent to the 
statement that $f'/f$ is constant. Bloch's principle thus suggests the
following normal families analogue proved by  
W.~Schwick~\cite[Theorem~5.1]{Schw89}.
\begin{thm}
{\bf Schwick's Theorem.}
Let $k\geq 2$ and let
${\mathcal F}$ be a family of functions holomorphic in 
a domain $D$. Suppose that $f$ and $f^{(k)}$ have
no zeros in $D$, for all $f\in{\mathcal F}$.
Then $\{f'/f:f\in\FF\}$ is normal.
\end{thm}
Theorem~\ref{thm11} was extended to meromorphic functions by
G.~Frank~\cite{Fra76}
for $k\geq 3$ and
by J.~K.~Langley~\cite{Lan93}
for $k=2$.
\begin{thm}
{\bf Frank-Langley Theorem.}
Let $f$ be meromorphic in $\C$ and let $k\geq 2$.
Suppose that $f$ and $f^{(k)}$ have no
zeros. Then $f$ has the form
$f(z)=e^{a z+b}$ or
$f(z)=(a z+b)^{-n},$
where $a,b\in\C$, $a\neq 0$, and $n\in\N$.
\end{thm}
We note that if $f$ has the form $f(z)=(a z+b)^{-n},$
then $f'(z)/f(z)=-n/(z+b/a)$ is not constant. 
So the original Bloch principle does not suggest that there is a normal family
analogue.
However, the family of all functions $f'/f$ of the above form is normal, and 
thus the following extension of Schwick's Theorem proved in~\cite{BL2}
is in accordance
with the variant of Bloch's principle discussed in~\S\ref{diffpol}.
\begin{thm} \label{thmbl}
{\bf Theorem.}
Let $k\geq 2$ and let
${\mathcal F}$ be a family of functions meromorphic in 
a domain $D$. Suppose that $f$ and $f^{(k)}$ have
no zeros in $D$, for all $f\in{\mathcal F}$.
Then $\{f'/f:f\in\FF\}$ is normal.
\end{thm}
For $k=2$ the result had been obtained already earlier in~\cite{Ber01}.
We omit the proofs of the above results and refer to the papers mentioned.

Instead of considering the condition that $f^{(k)}$ has no zeros one may,
more generally, consider the condition that 
$$L(f):=f^{(k)}+a_{k-1} f^{(k-1)}+\dots a_1 f'+a_0 f$$
has no zeros, for certain constants or functions $a_0,a_1,\dots,a_{k-1}$.
For functions 
meromorphic in the plane this has been addressed in~\cite{Bru92,Lan93,Lan94,FH,Ste87},
and results about normality appear in~\cite{Cli}.

Similarly one may replace the exceptional values of $f^{(k)}$ 
by exceptional values of $L(f)$
in many of the results discussed in \S\ref{intro2}-\S\ref{multiple}.
In fact, already in 1940 it was proved by C.~T.\ Chuang \cite{Chu40} that
in Miranda's Theorem~\ref{thmburmir}
one may replace the condition $f^{(k)}\neq 1$
by $L(f)\neq 1$ if the $a_j$ are holomorphic. We note that this
result can be deduced from the 
P\'olya-Saxer-Ullrich Theorem using 
the Zalcman-Pang Lemma in the same way  Miranda's Theorem 
was proved in~\S\ref{genzal}.
There are a large number of papers concerning 
exceptional values of $L(f)$. Here we only refer to~\cite{Dra69,Schi}.

\section{Fixed points and periodic points}\label{fixed}
\subsection{Introduction}\label{intro3}
Let $X,Y$ be sets,
let $f:X\to Y$ be a function,
and define
the iterates $f^n:X_n\to Y$ by $X_1:=X, f^1:=f$ and
$X_n:=f^{-1}(X_{n-1}\cap Y), f^n:=f^{n-1}\circ f$ for $n\in\N$, $n\geq 2$.
Note that $X_2=f^{-1}(X_1\cap Y)\subset X=X_1$ and thus $X_{n+1}\subset
X_n\subset X$ for all $n\in\N$.

A point $\xi\in X$ is called a {\em periodic point of period} $p$ of $f$
if $\xi\in X_p$ and $f^p(\xi)=\xi$, but $f^m(\xi)\neq \xi$ for $1\leq m\leq p-1$.
A periodic point of period 1 is called a {\em fixed point}.
The periodic points of period $p$ are thus the fixed points of $f^p$
which are not fixed points of $f^m$ for any $m$ less than $p$.
The periodic points play an important role in complex dynamics.

By the Fundamental Theorem of Algebra, every nonconstant 
polynomial $f$ which is not
of the form $f(z)=z+c$, $c\in\C\setminus\{0\}$
has a fixed point, and so does every iterate of $f$.
Transcendental entire functions need not have fixed points, as shown by
the example $f(z)=z+e^z$. P.\ Fatou~\cite[p.~345]{Fat26} proved
that the second iterate of a transcendental entire function 
has a fixed point, and this result was crucial for his proof that the 
Julia set of such a function is non-empty. P.\ C.\ Rosenbloom~\cite{Ros}
proved that in fact any iterate of a transcendental entire function
has infinitely fixed points. Summarizing these results we obtain the following 
theorem.
\begin{thm} \label{thmfatros}
{\bf Fatou-Rosenbloom Theorem.}
Let $f$ be a entire function and $p\in\N$, $p\geq 2$.
If $f^p$ has no fixed point, then $f$ has
the form $f(z)=z+c$ where $c\in\C\setminus\{0\}$.
\end{thm}
We note that the family of all functions $f$ of the form $f(z)=z+c$,
with $c\in\C\setminus\{0\}$,
is normal in $\C$. Thus the variant of Bloch's principle discussed in~\S\ref{diffpol} suggests
a normal family analogue not incorporated in
the original Bloch principle. This normal family analogue was proved by 
M.\ Ess\'en and S.\ Wu~\cite{Ess98},
thereby answering a question of  L.\ Yang~\cite[Problem 8]{Yan92}.
\begin{thm}\label{thmesswu}
{\bf Ess\'en-Wu Theorem.}
Let $D\subset \C$ be a domain and let $\FF$ be the family of all
holomorphic functions $f:D\to\C$ for which there exists $p=p(f)>1$ 
such that $f^p$ has no fixed point. Then $\FF$ is normal.
\end{thm}
We shall sketch the proof of the Ess\'en-Wu Theorem after Theorem~\ref{thmpropQ2}.

\subsection{Periodic points and quasinormality}\label{quasinorm}
In \S\ref{intro3} we discussed the property $P$ defined by
$\langle f,D\rangle \in P$ if $f$ is 
holomorphic in $D$ and if there exists $p=p(f)>1$ 
such that $f^p$ has no fixed point in $D$. 
The Ess\'en-Wu Theorem~\ref{thmesswu}
says that
$P$ implies normality.

We shall now be concerned with the weaker property $Q$ defined by
$\langle f,D\rangle \in Q$ if $f$ is 
holomorphic in $D$ and if there exists $p=p(f)>1$ 
such that $f$ has no periodic point of period $p$ in $D$.
We note that for $n\in\N$ the function
$f_n(z):=nz$ has no periodic points of period greater than $1$ 
so that $\langle f_n,\C\rangle \in Q$ for all $n$, but the $f_n$ do
not form a normal family. 

We first discuss what entire functions have property $Q$. 
For polynomials we have the following result due to I.\ N.\ Baker~\cite{Bak64}.
\begin{thm}{\bf Baker's Theorem.}
Let $f$ be a polynomial of degree $d\geq 2$ and let $p\in\N$, $p\geq 2$.
Suppose that  $f$ has no periodic point of period $p$.
Then $d=p=2$. Moreover, there exists a linear transformation $L$ 
such that $f(z)=L^{-1}(g(L(z)))$, with $g(z)=-z+z^2$.
\end{thm}
Note that for $g(z)=-z+z^2$ we have $g(z)-z=z(z-2)$ and 
$g^2(z)-z=z^3(z-2)$ so that there are no periodic points of period~$2$.

The case of transcendental functions is covered by the following generalization
of the Fatou-Rosenbloom Theorem~\ref{thmfatros} which was conjectured in
\cite[Problem 2.20]{Hay67} and proved in~\cite[Theorem~1]{Ber91}
and~\cite[\S 1.6, Satz~2]{Ber91a}.
\begin{thm}\label{thmber91}{\bf Theorem.}
Let $f$ be a transcendental entire function and let $p\in\N$, $p\geq 2$.
Then $f$ has infinitely many periodic points of period $p$.
\end{thm}
These results may be summarized as follows.
\begin{thm}{\bf Theorem.}\label{thmpropQ}
Let $\langle f,\C\rangle \in Q$. 
Then $f$ is a polynomial of degree at most $2$.
If $f$ has degree $2$, then
there exists a linear transformation $L$ 
such that $f(z)=L^{-1}(g(L(z)))$, with $g(z)=-z+z^2$.
\end{thm}
As mentioned above, property $Q$ does not imply normality.
However, we have the following result~\cite{BB}.

\begin{thm}{\bf Theorem.}\label{thmpropQ2}
For every
domain $D\subset\C$ the family
 $\{f:\langle f,D\rangle \in Q \}$ is quasinormal of order~$1$ in $D$. 
\end{thm}

The proofs of Theorem~\ref{thmpropQ2} and the Ess\'en-Wu Theorem~\ref{thmesswu} are
based on similar arguments. 

{\em We sketch the idea.} For simplicity we only prove that $\FF$ is normal
if $f^2$ has no fixed point for all $f\in\FF$, and that $\FF$ is quasinormal
of order $1$ if $f$ has no periodic point of period $2$ for all $f\in\FF$. 
The general case is proved along the same lines; we refer to the papers cited
for the details.

First we prove that a family $\FF$ of functions holomorphic in a domain $D$
is quasinormal
of order $3$ if $f$ has no periodic point of period $2$ for all $f\in\FF$.
Suppose that $\FF$ is not quasinormal of order $3$.
Then there exists a sequence $(f_n)$ in $\FF$ and 
four points 
$a_1,a_2,a_3,a_4\in D$ such that no subsequence of $(f_n)$ is normal at any of
the points $a_j$. 

Applying Ahlfors's Theorem~\ref{ahlspec} with a domain $D_3$ containing $\infty$ 
we see that if $D_1,D_2$ are Jordan domains in $\C$  with disjoint closures,
if $\Omega$ is a neighborhood of one of the points $a_j$, and if $n$ is sufficiently large,
then $f_n$ has an island $U$  contained in $\Omega$ over one of the domains $D_1$ or $D_2$.
We choose $\varepsilon >0$ such that the closures of the disks of radius $\varepsilon $
around the $a_j$ are pairwise disjoint.  We see that if $n$ is sufficiently large and
$j,k_1,k_2\in\{1,2,3,4\}$ with $k_1\neq k_2$, 
then $f_n$ has an island $U$ in $D(a_j,\varepsilon)$ over
$D(a_{k_1},\varepsilon)$ or $D(a_{k_2},\varepsilon)$. Thus $f_n$ has an island in $D(a_j,\varepsilon)$ 
over $D(a_{k},\varepsilon)$ for at least three values of $k$.
This implies that there exists $j,k\in\{1,2,3,4\}$, $j\neq k$, such that 
$f_n$ has an island $U$ in $D(a_j,\varepsilon)$ over $D(a_{k},\varepsilon)$ 
and an island $V$ in $D(a_k,\varepsilon)$ over $D(a_{j},\varepsilon)$.
We now consider a component $W$ of $U\cap f_n^{-1}(V)$ and see that
$f_n^2|_W:W\to D(a_{j},\varepsilon)$ is a proper map.
In particular, $f_n^2$ takes the value $a_j$ in $W$.

For $z\in\partial W$ we have 
$$|(f_n^2(z)-a_j)-\left(f_n^2(z)-z\right)|=
|z-a_j|<\varepsilon=|f_n^2(z)-a_j|.$$
Rouch\'e's Theorem implies that $f_n^2(z)-z$ has a zero in $W$,
say $f^2(\xi)=\xi$ where $\xi\in W$. Since $f(\xi)\in D(a_k,\varepsilon )$
and $W\cap D(a_k,\varepsilon )=\emptyset$ we see that $\xi$ is a periodic 
point of period~$2$. This contradicts our assumption. Thus 
$\FF$ is quasinormal
of order~$3$ if $f$ has no periodic point of period $2$ for all $f\in\FF$ --
and thus in particular if  $f^2$ has no fixed point for all $f\in\FF$. 

To complete the proof that $\FF$ is normal if $f^2$ has no fixed point for all 
$f\in\FF$, suppose that $\FF$ is not normal. Then there exists a sequence
$(f_n)$ in $\FF$ and a point 
$a_1\in D$ such that no subsequence of~$(f_n)$ is normal at~$a_1$. 
Since $\FF$ is quasinormal of order~$3$ we may, passing to a subsequence
if necessary, assume that $f_n$ converges in $D\setminus\{a_1,a_2,a_3\}$
where $a_2,a_3\in D$.
The maximum principle implies that $f_n\to\infty$ in 
$D\setminus\{a_1,a_2,a_3\}$.
We find that if $\varepsilon >0$ is such that the closure of the disk 
$D(a_1,\varepsilon )$ is contained in $D\setminus\{a_2,a_3\}$ 
and if $n$ is large enough,
then $f_n$ has an island $U$ in $D(a_1,\varepsilon )$ 
over $D(a_1,\varepsilon )$.
As above Rouch\'e's Theorem implies that $f_n$ has a fixed point in $U$.
This fixed point is also a fixed point of $f_n^2$, contradicting the assumption.

The proof  that $\FF$ is quasinormal
of order $1$ if $f$ has no periodic point of period $2$ for all $f\in\FF$
is completed in a similar fashion. Assuming that this is not the case we find
a sequence
$(f_n)$ in $\FF$ and 
two points 
$a_1,a_2\in D$ such that no subsequence of $(f_n)$ is normal at $a_1$ or $a_2$.
Passing to a subsequence we may again assume that 
$f_n$ converges in $D\setminus\{a_1,a_2,a_3\}$ for some $a_3\in D$,
and hence
$f_n\to\infty$ in $D\setminus\{a_1,a_2,a_3\}$ by the maximum principle.
For suitable $\varepsilon >0$ and sufficiently large $n$ we find that
$f_n$ has an island $U$ in $D(a_1,\varepsilon)$ over $D(a_{2},\varepsilon)$ 
and an island $V$ in $D(a_2,\varepsilon)$ over $D(a_{1},\varepsilon)$.
Again we consider a component $W$ of $U\cap f_n^{-1}(V)$ and see that
$f_n^2|_W:W\to D(a_{1},\varepsilon)$ is a proper map.
As above, we see that $W$ contains a periodic point $\xi$ of period $2$ of~$f_n$.\qed

A fixed point $\xi$ of a holomorphic function $f$ is called {\em repelling}
if $|f'(\xi)|>1$. Repelling periodic points are defined accordingly.
They play an important role in complex dynamics.
Many of the results mentioned above have generalizations where instead 
of fixed points and periodic points only repelling fixed points and periodic points 
are considered. For example, the Theorems~\ref{thmesswu} and~\ref{thmber91} hold
literally with the word ``repelling'' added. But the results about polynomials 
are somewhat different; see~\cite{Ber91,Ber91a,Berquasi,Ess} for more details.

The condition that $f$ has no (repelling)
periodic points of some period --
or that some iterate does not have (repelling) fixed points -- has also been considered 
for meromorphic functions. We refer 
to~\cite{CF05,Ess98,Sie,Sie2,WangWu} for results when this implies
(quasi)normality,  to~\cite{Bak64,Kis} for results concerning rational functions,
and to~\cite[\S 3]{Ber93} 
for the case of transcendental functions meromorphic in the plane.
\section{Further topics}\subsection{Functions sharing values}\label{sharedvalues}
Two meromorphic functions $f$ and $g$  are said to share a value $a\in\CC$ if
they have the same $a$-points; that is,
$f(z)=a$ if $g(z)=a$ and vice versa.
A famous  result of Nevanlinna \cite{Nev26}
says that if two functions meromorphic in the
plane share five values, then they are equal.
E.~Mues and N. Steinmetz~\cite{MueSte} proved that if $f$ is meromorphic in the plane
and if $f$ and $f'$ share three values, then $f'=f$ so that $f(z)=c e^z$ for some 
$c\in\C$. Now the family $\{c e^z:c\in\C\}$  is normal and thus the variant of 
Bloch's principle discussed in~\S\ref{diffpol} suggests that the family of functions $f$
meromorphic in a domain and sharing three fixed values with their derivative
is normal. W.~Schwick~\cite{Schw92} proved that this is in fact the case.

On the other hand, G.~Frank and W.~Schwick~\cite{FS1,FS2} showed that 
for a function $f$ meromorphic in the plane 
the condition that $f$ and $f^{(k)}$ share three values for some
$k\geq 2$ still implies that $f=f^{(k)}$, but the family of all
functions $f$  which are meromorphic in some domain and satisfy this
condition is not normal.
This is in accordance with 
with both the original Bloch principle and its 
variant 
introduced in~\S\ref{diffpol}, since the 
functions $f$ which are meromorphic in the plane and satisfy $f=f^{(k)}$
do not form a normal family for~$k\geq 2$.

There is an enormous amount of
literature on functions meromorphic in the plane that share values,
and in recent years many papers on corresponding normality results have appeared.
Here we only refer to~\cite{FZ,Pang02,PZ,PZ2} and the literature cited there.
We note that some of these results generalize the results about exceptional
values of derivatives described in~\S\ref{intro2}, since if two functions
omit the same value, then they of course also share this value.
\subsection{Gap series}\label{gaps}
A classical result of L.~Fej\'er \cite[p.~412]{Fej08} says that an entire
function $f$ of the form 
\begin{equation}\label{fejer}
f(z)=\sum_{k=0}^\infty a_k z^{n_k}
\ \ \ \text{where} \ \ \ 
\sum_{k=0}^\infty \frac{1}{n_k}=\infty
\end{equation}
has at least one zero. 
St.\ Ruscheweyh and K.-J.\ Wirths~\cite{RuWi} have shown 
that the family of all functions $f$ of the form (\ref{fejer}) which
are holomorphic in the unit disk and do not vanish there
form a normal family. 

There are a number of further results, as well as open questions,
 on exceptional values of entire
functions with gap series; see~\cite{Mur} for a survey. 
Here we only mention a  question of G.\ P\'olya~\cite[p.~639]{Pol} whether 
Fej\'er's condition $\sum_{k=0}^\infty 1/n_k=\infty$ in~(\ref{fejer})
can be replaced by Fabry's condition $\lim_{k\to\infty} n_k/k=\infty$.

In accordance with Bloch's principle, 
Ruscheweyh and Wirths~\cite{RuWi} have raised the following conjecture.
\begin{exam}
{\bf Conjecture.}
For $\Lambda\subset\N\cup\{0\}$ and a 
function $f$ holomorphic in a domain $D$ containing $0$ define 
$\langle f,D\rangle\in P_\Lambda$ if $f$ has a power series expansion
$$f(z)=\sum_{\lambda\in\Lambda}a_\lambda z^\lambda$$
and if $f(z)\neq 0$ for all $z\in D$.
Then $\{f:\langle f,\D\rangle\in P_\Lambda\}$ is normal in $\D$ if and only if 
$\langle f,\C\rangle\in P_\Lambda$ implies that $f$ is constant.
\end{exam}
We note that one direction in this conjecture is obvious:
Normality of the family $\{f:\langle f,\D\rangle\in P_\Lambda\}$ implies that 
$\langle f,\C\rangle\in P_\Lambda$ only for constant functions $f$.
In fact, if a nonconstant 
entire function $f$ has property  $P_\Lambda$, then so does every function in
the family $\{f(nz):n\in\N\}$, and this family is not normal at $0$.

We mention that Zalcman's Principle~\ref{zalcprinc}
cannot apply to properties concerning gap series 
since condition $(ii)$ is not satisfied.
However, there are some further results in addition to~\cite{RuWi}
which support the above conjecture.
W.~K.\ Hayman~\cite{Hay83} has considered entire functions
with gaps in arithmetic progressions. Normal family analogues of 
some of the results have been obtained St.~Ruscheweyh and L.~Salinas~\cite{RuSa}
and by J.\ Grahl~\cite{Gra00}.
\subsection{Holomorphic curves}\label{holocurves}
It is easily seen that Picard's Theorem is equivalent to 
the statement that if $f_1,f_2,f_3$ are nonvanishing entire functions 
and if 
$c_1,c_2,c_3$ are nonzero complex numbers
such that $\sum_{j=1}^3 c_j f_j=0$,
then each quotient $f_j/f_k$ is constant.
In fact, 
writing $F:=-c_1f_1/(c_3f_3)$
we see that $F$ is an entire function without zeros. Moreover, since
$F=1+c_2f_2/(c_3f_3)$ we see that $F$ also omits the value $1$. Thus
$F$ is constant by Picard's Theorem.
This implies that not only $f_1/f_3$ but also the other quotients $f_1/f_2$ and 
$f_2/f_3$ are constant.

Another way to phrase this result is that the hypothesis that the three functions
$f_1,f_2,f_3$ are linearly dependent already implies that two of them are
linearly dependent.

A generalization of this statement was proved by \'E.~Borel~\cite{Borel97}.
\begin{thm}
{\bf Borel's Theorem.}
Let $p\in\N$, let $f_1,\dots,f_p$ be  entire functions 
without zeros 
and let $c_1,\dots,c_p\in\C\setminus\{0\}$.
Suppose that $\sum_{j=1}^p c_j f_j=0$.
Then $\{f_1,\dots,f_p\}$ contains a linearly dependent subset of less than
$p$ elements.
\end{thm}
Repeated application of this result shows that under the hypotheses of Borel's Theorem
the set $\{1,\dots,p\}$ can be written as the union of disjoint subsets $I_\mu$,
each of which has at least two elements, and such
that if $j,k$ are in the same set $I_\mu$, then $f_j/f_k$ is constant.

One may ask whether this Picard type theorem also has an analogue in the context of normal families.
However, as noted already by  Bloch~\cite[p.~311]{Blo26b} himself,
it is not clear at first sight what such an analogue could look like.
This problem
was then addressed by H.~Cartan~\cite{Cartan28}
who proved such an analogue in the case that $p=4$ and made a conjecture for the general 
case. Cartan's conjecture was disproved by A.~Eremenko~\cite{Ere96a}. However, Eremenko~\cite{Ere96b}
also showed that a weakened form of the conjecture is true for $p=5$. For a connection to
gap series we refer to a paper by
J.~Grahl~\cite{Gra00}.

\subsection{Quasiregular maps}\label{quasireg}
Many of the concepts used in this survey are
also applicable for quasiregular maps in higher dimensions; see
\cite{Rick} for the definition and basic properties of quasiregular maps.
One of the most important results in this theory is the analogue of
Picard's Theorem which was obtained by S.~Rickman~\cite{Ric80}.
He proved
that there exists $q=q(d,K)\in\N$ with the property that every $K$-quasiregular
map $f:\R^d\to\R^d$ which omits $q$ points is constant.
The corresponding normality result was proved by R.~Miniowitz~\cite{Min}, using 
an extension of the Zalcman Lemma to quasiregular maps.
A.~Eremenko~\cite{Ere00}  has used Miniowitz's result 
to extend the classical covering theorem of Bloch~\cite{Blo26c}
to quasiregular maps.

Recall that the limit function $f$ occuring in Zalcman's lemma
has bounded spherical derivative. 
The corresponding conclusion in the context of quasiregular maps
is that the limit function is uniformly continuous.
This plays an important role in the work 
of M.~Bonk and J.~Heinonen~\cite{BH} 
on closed, connected and oriented Riemannian $d$-manifolds~$N$ 
for which there exists a nonconstant quasiregular map $f:\R^d\to N$. 
They use the above argument to show that if there is such a 
mapping, then there is also a uniformly continuous one.

Many of the results of~\S\ref{intro3} and \S\ref{quasinorm} concerning fixed points have 
been extended to quasiregular mappings by H.~Siebert~\cite{Sie,Sie2}.
For example, Theorems~\ref{thmesswu}, \ref{thmber91} and~\ref{thmpropQ2} hold literally
for quasiregular maps. We omit a detailed discussion of the
results about quasiregular mappings  here and refer to the papers cited.

\begin{acknowledgement}
This paper is based on the Lecture Notes for the Advanced Course in Operator
Theory and Complex Analysis in  Sevilla in June 2005.
I~would like to thank Alfonso Montes Rodr{\i}guez for the invitation to this interesting 
workshop. I am also indebted to David Drasin, Alex Eremenko,
Jim Langley and Larry Zalcman
for some helpful comments on these notes. 
\end{acknowledgement}

\noindent 
\end{document}